\documentclass[11pt,a4paper,twoside]{article}
\usepackage{amssymb}
\usepackage{amsmath}
\usepackage{fancyhdr}
\usepackage{amsthm}
\usepackage{rotate}
\usepackage{graphicx}
\usepackage[T1]{fontenc}

\usepackage{afterpage}  %

\newtheorem{theorem}{Theorem}[section]
\theoremstyle{plain}

\newtheorem{remark}[theorem]{Remark}

\numberwithin{equation}{section}

\newenvironment{GAPexample}
{\vskip 3mm\begin{scriptsize}} {\end{scriptsize}\vskip 1mm}

\def\returns{\ $\triangleright$\ }         

\def\myline#1{\newline\phantom{xx}#1}
\def\myextraline#1{\newline\phantom{xxxx}#1}

\begin{document}

\afterpage{\rhead[]{\thepage} \chead[\small G. P. Nagy and P. Vojt\v{e}chovsk\'y         
]{\small  Computing with loops} \lhead[\thepage]{} }                  

\begin{center}
\vspace*{2pt}
{\Large \textbf{Computing with small quasigroups and loops}}\\[36pt]
{\large \textsf{\emph{G\'abor P.~Nagy and Petr Vojt\v{e}chovsk\'y}}}
\\[36pt]
\textbf{Abstract}
\end{center}
{\footnotesize This is a companion to our lectures \emph{GAP and loops}, to be
delivered at the \emph{Workshops Loops 2007}, Prague, Czech Republic. In the
lectures we introduce the GAP \cite{GAP} package LOOPS \cite{LOOPS}, describe
its capabilities, and explain in detail how to use it. In this paper we first
outline the philosophy behind the package and its main features, and then we
focus on three particular computational problems: construction of loop
isomorphisms, classification of small Frattini Moufang loops of order $64$, and
the search for loops of nilpotency class higher than two with an abelian inner
mapping group.

In particular, this is not a user's manual for LOOPS, which can be downloaded
from the distribution website of LOOPS.}

\footnote{\textsf{2000 Mathematics Subject Classification:} Primary 20N05.}

\footnote{\textsf{Keywords:} loop, quasigroup, GAP, computation
in nonassociative algebra, Latin square, loop isomorphism, Cs\"{o}rg\H{o} loop,
small Frattini Moufang loop, code loop, LOOPS package.}

\footnote{This paper was written during the Marie Curie Fellowship of the first
author at the University of W\"{u}rzburg. The second author supported by the
PROF 2006 grant of the University of Denver.}

\section{Main features}

On the one hand, since there is no useful representation theory for quasigroups
and loops, we have decided to represent quasigroups and loops in LOOPS by their
Cayley tables, thus effectively limiting the scope of the package to
quasigroups of order at most 500 or so. (A future project is to implement other
loop representations, notably by connected transversals in groups.)

On the other hand, to take advantage of the powerful methods for groups already
present in GAP, most calculations in LOOPS are delegated to the permutation
groups associated with quasigroups, rather than performed on the level of
Cayley tables. For instance, to decide if a loop is simple, we check whether
its multiplication group is a primitive permutation group.

To avoid repeated calculations, we store most information obtained for a given
quasigroup as its attribute. In GAP, there is no syntactical difference between
calling a method or retrieving an attribute. For instance, when $Q$ is a
quasigroup then \texttt{Center($Q$)} calculates and stores the center $Z(Q)$ of
$Q$ when called for the first time, while it retrieves the stored attribute
$Z(Q)$ when called anytime later.

Moreover, GAP uses simple deduction process---filters---to obtain additional
information about an object without an explicit user's request. For instance,
if LOOPS knows that $Q$ is a left Bol loop that is also commutative, the
built-in filter \texttt{(IsMoufangLoop, IsLeftBolLoop and IsCommutative)}
automatically deduces that $Q$ is a Moufang loop and stores this information
for $Q$. This is a powerful tool, since many filters built into LOOPS are deep
theorems.

\subsection{Creating quasigroups and loops} A \emph{$($quasigroup$)$ Cayley table} is
an $n\times n$ Latin square with integral entries $x_1<\cdots <x_n$. A
\emph{canonical Cayley table} is a Cayley table with $x_1=1$, $\dots$, $x_n=n$.

When $T$ is a Cayley table, \texttt{QuasigroupByCayleyTable($T$)} creates a
quasigroup whose Cayley table is the canonical Cayley table obtained from $T$
by replacing $x_i$ with $i$. Should $T$ be \emph{normalized}---the first row
and first column reads $x_1$, $\dots$, $x_n$---then
\texttt{LoopByCayleyTable($T$)} returns the corresponding loop. The Cayley
table of a quasigroup $Q$ can be retrieved by \texttt{CayleyTable($Q$)}.

Throughout this paper, we illustrate the methods of LOOPS by examples, often
without any comments for self-explanatory commands. The syntax is that of GAP.

\begin{GAPexample}
\begin{verbatim}
gap> Q := QuasigroupByCayleyTable([[2,1],[1,2]]); Elements(Q);
<quasigroup of order 2>
[ q1, q2 ]
gap> L := LoopByCayleyTable([[3,5],[5,3]]); Elements(L); L.2;
<loop of order 2>
[ l1, l2 ]
l2
gap> CayleyTable(Q);
[ [ 2, 1 ], [ 1, 2 ] ]
gap> Print(L);
<loop with multiplication table
[ [  1,  2 ],
  [  2,  1 ] ]
>
\end{verbatim}
\end{GAPexample}

It is also possible to create quasigroups and loops by reading Cayley tables
from files (with very relaxed conditions on the form of the Cayley table), by
converting groups to quasigroups, by taking subquasigroups, subloops, factor
loops, direct products, etc. See the manual for details.

\subsection{Conversions}

Even if a quasigroup happens to have a neutral element, it is not considered a
loop in LOOPS unless it is declared as a loop. Similarly, a group of GAP is not
considered a loop. We therefore provide conversions between these types of
algebras:

\begin{GAPexample}\begin{verbatim}
gap> G := Group((1,2,3),(1,2)); AsLoop(G);
Group([ (1,2,3), (1,2) ])
<loop of order 6>
gap> Q := QuasigroupByCayleyTable([[2,1],[1,2]]);  AsLoop(Q);
<quasigroup of order 2>
<loop of order 2>
\end{verbatim}\end{GAPexample}

The neutral element of any loop $L$ in LOOPS is always the first element of
$L$, i.e., \texttt{One($L$) = $L$.1}.

Given a quasigroup $Q$ and elements $f$, $g\in Q$, the principal loop isotope
$(Q,f,g)$ of $Q$ is obtained from $Q$ via the isotopism $(R_g^{-1}, L_f^{-1},
\textrm{id})$, cf. \cite[p. 60]{Pflugfelder}. Then $(Q,f,g)$ is a loop with
neutral element $fg$.

The conversion \texttt{AsLoop($Q$)} works as follows, starting with a
quasigroup $Q$:
\begin{enumerate}
\item[(i)] When $Q$ does not have a neutral element, it is first replaced by the
principal loop isotope $(Q,Q.1,Q.1)$, thus turning $Q$ into a loop with neutral
element $(Q.1)(Q.1)$.
\item[(ii)] When $Q$ has a neutral element $k$, it is replaced by its isomorphic
copy via the transposition $(1,k)$.
\end{enumerate}

\subsection{Subquasigroups and subloops}

A new quasigroup $Q_2$ is frequently obtained as a subquasigroup of an existing
quasigroup $Q_1$. Since all information about $Q_2$ is already contained in the
Cayley table of $Q_1$, and since it is often desirable to have access to the
embedding of $Q_2$ into $Q_1$, we provide a mechanism in LOOPS for maintaining
the inclusion of $Q_2$ and $Q_1$.

When $Q_1$ is a quasigroup and $S$ is a subset of $Q_1$,
\texttt{Subquasigroup($Q_1$, $S$)} returns the subquasigroup $Q_2$ of $Q_1$
generated by $S$. At the same time, the attribute \texttt{Parent($Q_2$)} is set
to \texttt{Parent($Q_1$)}, hence ultimately pointing to the largest quasigroup
from which $Q_2$ has been created. The elements of $Q_2$ and the Cayley table
of $Q_2$ are then calculated relative to the parent of $Q_2$.

\begin{GAPexample}\begin{verbatim}
gap> L := AsLoop(Group((1,2,3),(1,2))); S := Subloop(L,[3]);
<loop of order 6>
<loop of order 2>
gap> Parent( S ) = L; PosInParent( S ); Elements( S );
true
[ 1, 3 ]
[ l1, l3 ]
gap> HasCayleyTable( S ); CayleyTable( S );
false
[ [ 1, 3 ], [ 3, 1 ] ]
\end{verbatim}\end{GAPexample}

Note that the Cayley table of a subquasigroup is created only upon user's request.

\subsection{Bijections as permutations on $\{1,\dots,n\}$}

When calculating isomorphisms, isotopisms, or other bijections of quasigroups
of order $n$, the result is always returned as a permutation (triple of
permutations) of $\{1,\dots,n\}$. Equivalently, the quasigroups in question are
first replaced by isomorphic copies with canonical Cayley tables, and only then
the bijections are calculated. It is always possible to reconstruct the
original bijection using the attribute \texttt{PosInParent}.

\subsection{A few words about the implementation}

One of the biggest strengths of the computer algebra system GAP is that most
algebraic structures can be defined within it. In this subsection we briefly
explain how the variety of quasigroups is implemented in LOOPS. In order to
understand the implementation, we will need the following GAP terminology:

\begin{itemize}
\item[-] A \textit{filter}, such as \texttt{IsInteger} and
\texttt{IsPermGroup}, is a special unary function on the set of GAP objects
which returns either \texttt{true} or \texttt{false}. Roughly speaking, a
filter is an \textit{a priori} attribute of an object.

\item[-] A \textit{category} is a class of objects defined by a collection of
filters. An object can lie in several categories. For example, a row vector
lies in the categories \texttt{IsList} and \texttt{IsVector}.

\item[-] All GAP objects are partitioned into \textit{families.} The family of
an object determines its relation to other objects. For instance, all
permutations form a family, and groups presented by generators and relations
form another family. However, a family is not a collection of objects, but
abstract information about objects.

\item[-] Beside its name, a family can have further \textit{labels}.

\item[-] Every GAP object has a \textit{type.} The type of an object determines
if a given operation can be performed with that object, and if so, how it is to
be performed. The type of an object is derived from its family and its filters.

\item[-] A given data structure can be made into an \textit{object} by
specifying its type, that is, its family and its filters.
\end{itemize}

The following function constructs a quasigroup $Q$ with Cayley table
\texttt{ct}. First we define a family corresponding to the elements of $Q$ and
tell GAP that it will consist of quasigroup elements. Then we objectify the
individual elements in this family, and label the family by the set of its
elements, by the size of $Q$, and by the Cayley table. Then we objectify $Q$
whose family will be the \textit{collection} of its elements. Finally, we set
some important attributes of $Q$.

\begin{GAPexample}\begin{verbatim}
function( ct )
    local F, Q, elms, n;
    # constructing the family of the elements of this quasigroup
    F := NewFamily( "QuasigroupByCayleyTableFam", IsQuasigroupElement );
    # installing data ("labels") for the family
    n := Length ( ct );
    F!.size := n;
    elms := Immutable( List( [1..n], i -> Objectify(
        NewType( F, IsQuasigroupElement and IsQuasigroupElmRep), [ i ] ) ) );
    F!.set := elms;
    F!.cayleyTable := ct;
    # creating the quasigroup by turning it into a GAP object
    # the family of Q is the collection of its elements
    Q := Objectify( NewType( FamilyObj( elms ),
        IsQuasigroup and IsAttributeStoringRep ), rec() );
    # setting some attributes for the quasigroup
    SetSize( Q, n );
    SetAsSSortedList( Q, elms );
    SetCayleyTable( Q, ct );
    return Q;
end;
\end{verbatim}\end{GAPexample}

Operations in GAP are overloaded, i.e., the same operation can be applied to
different types of objects. In order to deal with this situation, GAP uses a
method selection: When an operation is called, GAP first checks the types of
the arguments, and then selects the appropriate method.

Here is how the multiplication of two quasigroup elements is implemented:

\begin{GAPexample}\begin{verbatim}
InstallMethod( \*, "for two quasigroup elements",
    IsIdenticalObj,
    [ IsQuasigroupElement, IsQuasigroupElement ],
function( x, y )
    local F;
    F := FamilyObj( x );
    return F!.set[ F!.cayleyTable[ x![ 1 ] ][ y![ 1 ] ] ];
end );
\end{verbatim}\end{GAPexample}

Note that the underlying quasigroup is easily accessed since the element $x$
\emph{knows} into which quasigroup it belongs.

\section{What is in the package}

Here is a very brief overview of the methods implemented in LOOPS, version
1.4.0. See the manual for (much) more details. Argument $Q$ stands for a
quasigroup, and $L$ for a loop. Thus the methods with argument $Q$ apply to
both quasigroups and loops, while those with argument $L$ apply only to loops.
Any additional restrictions on the arguments are listed in parentheses. The
symbol\returns is a shortcut for \emph{returns}.

\subsection{Basic methods and attributes}

\noindent Cayley tables and elements:
\myline{\texttt{Elements($Q$)}\returns list of elements of $Q$,}
\myline{\texttt{CayleyTable($Q$)}\returns Cayley table of $Q$,}
\myline{\texttt{One($L$)}\returns the neutral element of $L$,}
\myline{\texttt{MultiplicativeNeutralElement($Q$)}\returns the neutral element of $Q$, or fail}
\myline{\texttt{Size($Q$)}\returns the size of $Q$,}
\myline{\texttt{Exponent($L$)}\returns the exponent of $L$ ($L$ power-associative).}

\noindent Arithmetic operations:
\myline{\texttt{LeftDivision($x,y$)}\returns $x{\setminus}y$,}
\myline{\texttt{RightDivision($x,y$)}\returns $x/y$,}
\myline{\texttt{LeftDivisionCayleyTable($Q$)}\returns Cayley table of left division in $Q$,}
\myline{\texttt{RightDivisionCayleyTable($Q$)}\returns Cayley table of right division in $Q$.}

\noindent Powers and inverses:
\myline{\texttt{LeftInverse($x$)}\returns $x^\lambda$, where $x^\lambda x=1$,}
\myline{\texttt{RightInverse($x$)}\returns $x^\rho$, where $xx^\rho=1$,}
\myline{\texttt{Inverse($x$)}\returns the two-sided inverse of $x$, if it exists.}

\noindent Associators and commutators:
\myline{\texttt{Associator($x,y,z$)}\returns the unique element $u$ with $(xy)z = (x(yz))u$},
\myline{\texttt{Commutator($x,y$)}\returns the unique element $v$ with $xy = (yx)v$.}

\noindent Generators:
\myline{\texttt{GeneratorsOfQuasigroup($Q$)}\returns a generating subset of $Q$,}
\myline{\texttt{GeneratorsOfLoop($L$)}\returns a generating subset of $L$,}
\myline{\texttt{GeneratorsSmallest(Q)}\returns a generating subset of $Q$ of size $\le\log_2|Q|$.}

\noindent Subquasigroups:
\myline{\texttt{IsSubquasigroup($Q,S$)}\returns true if $S$ is a subquasigroup of $Q$,}
\myline{\texttt{IsSubloop($L,S$)}\returns true if $S$ is a subloop of $L$,}
\myline{\texttt{AllSubloops($L$)}\returns list of all subloops of $L$,}
\myline{\texttt{RightCosets($L,S$)}\returns right cosets modulo $S$ ($S\le L$),}
\myline{\texttt{RightTransversal($L,S$)}\returns a right transversal modulo $S$ ($S\le L$)}.

\noindent Translations and sections:
\myline{\texttt{LeftTranslation($Q,x$)}\returns the left translation $L_x$ by $x$ in $Q$ ($x\in Q$),}
\myline{\texttt{RightTranslation($Q,x$)}\returns the right translation $R_x$ by $x$ in $Q$ ($x\in Q$),}
\myline{\texttt{LeftSection($Q$)}\returns the set of all left translations in $Q$,}
\myline{\texttt{RightSection($Q$)}\returns the set of all right translations in $Q$.}

\noindent Multiplication groups:
\myline{\texttt{LeftMultiplicationGroup($Q$)}\returns the left multiplication group of $Q$,}
\myline{\texttt{RightMultiplicationGroup($Q$)}\returns the right multiplication group of $Q$,}
\myline{\texttt{MultiplicationGroup($Q$)}\returns the multiplication group of $Q$,}
\myline{\texttt{RelativeLeftMultiplicationGroup($L,S$)}\returns the group generated by all}
\myextraline{left translations of $L$ restricted to $S$ ($S\le L$),}
\myline{\texttt{RelativeRightMultiplicationGroup($L,S$)}\returns the group generated by all}
\myextraline{right translations of $L$ restricted to $S$ ($S\le L$),}
\myline{\texttt{RelativeMultiplicationGroup($L,S$)}\returns the group generated by all}
\myextraline{translations of $L$ restricted to $S$ ($S\le L$).}

\noindent Inner mapping groups:
\myline{\texttt{InnerMappingGroup($L$)}\returns the inner mapping group of $L$,}
\myline{\texttt{LeftInnerMappingGroup($L$)}\returns the group generated by $L_{yx}^{-1}L_yL_x$,}
\myline{\texttt{RightInnerMappingGroup($L$)}\returns the group generated by $R_{xy}^{-1}R_yR_x$.}

\noindent Nuclei:
\myline{\texttt{LeftNucleus($Q$)}\returns the left nucleus of $Q$,}
\myline{\texttt{RightNucleus($Q$)}\returns the right nucleus of $Q$,}
\myline{\texttt{MiddleNucleus($Q$)}\returns the middle nucleus of $Q$,}
\myline{\texttt{Nuc($Q$)}, \texttt{NucleusOfQuasigroup($Q$)}\returns the nucleus of $Q$.}

\noindent Commutant, center and associator subloop:
\myline{\texttt {Commutant($Q$)}\returns $\{x\in Q;\;xy=yx$ for every $y\in Q\}$,}
\myline{\texttt{Center($Q$)}\returns the center of $Q$,}
\myline{\texttt{AssociatorSubloop($L$)}\returns the smallest $S\unlhd L$ such that $L/S$ is a group.}

\noindent Normal subloops:
\myline{\texttt{IsNormal($L,S$)}\returns true if $S$ is a normal subloop of $L$,}
\myline{\texttt{NormalClosure($L,S$)}\returns the smallest normal subloop of $L$ containing $S$,}
\myline{\texttt{IsSimple($L$)}\returns true if $L$ is a simple loop.}

\noindent Factor loops:
\myline{\texttt{FactorLoop($L,N$)}\returns $L/N$ ($N$ normal subloop of $L$),}
\myline{\texttt{NaturalHomomorphismByNormalSubloop($L,N$)}\returns the projection}
\myextraline{$L\to L/N$ ($N$ normal subloop of $L$).}

\noindent Central nilpotency and central series:
\myline{\texttt{NilpotencyClassOfLoop($L$)}\returns the (central) nilpotency class of $L$,}
\myline{\texttt{IsNilpotent($L$)}\returns true if $L$ is nilpotent,}
\myline{\texttt{IsStronglyNilpotent($L$)}\returns true if the mult. group of $L$ is nilpotent,}
\myline{\texttt{UpperCentralSeries($L$)}\returns the upper central series of $L$,}
\myline{\texttt{LowerCentralSeries($L$)}\returns the lower central series of $L$,}

\noindent Solvability:
\myline{\texttt{IsSolvable($L$)}\returns true if $L$ is solvable,}
\myline{\texttt{DerivedSubloop($L$)}\returns the derived subloop of $L$,}
\myline{\texttt{DerivedLength($L$)}\returns the derived length of $L$,}
\myline{\texttt{FrattiniSubloop($L$)}\returns the Frattini subloop of $L$ ($L$ strongly nilpotent).}

\noindent Isomorphisms and automorphisms:
\myline{\texttt{IsomorphismLoops($L,M$)}\returns an isomorphism of loops $L\to M$, or fail,}
\myline{\texttt{LoopsUpToIsomorphism($ls$)}\returns filtered list $ls$ of loops up to isomorphism,}
\myline{\texttt{AutomorphismGroup($L$)}\returns the automorphism group of $L$,}
\myline{\texttt{IsomorphicCopyByPerm($Q,p$)}\returns an isomorphic copy of $Q$ via the}
\myextraline{permutation $p$,}
\myline{\texttt{IsomorphicCopyByNormalSubloop($L,S$)}\returns an isomorphic copy of $L$ in}
\myextraline{which $S\unlhd L$ occupies the first $|S|$ elements of $L$ and where the remaining}
\myextraline{elements correspond to the cosets of $S$ in $L$.}

\noindent Isotopisms:
\myline{\texttt{IsotopismLoops($L,M$)}\returns an isotopism $L\to M$, or fail,}
\myline{\texttt{LoopsUpToIsotopism($ls$)}\returns filtered list $ls$ of loops up to isotopism.}

\subsection{Testing properties of quasigroups and loops}

Associativity, commutativity and generalizations:
\myline{\texttt{IsAssociative($Q$)}\returns true if $Q$ is associative,}
\myline{\texttt{IsCommutative($Q$)}\returns true if $Q$ is commutative,}
\myline{\texttt{IsPowerAssociative($L$)}\returns true if $L$ is power associative,}
\myline{\texttt{IsDiassociative($L$)}\returns true if $L$ is diassociative.}

\noindent Inverse properties:
\myline{\texttt{HasLeftInverseProperty($L$)}\returns true if $x^\lambda(xy)=y$,}
\myline{\texttt{HasRightInverseProperty($L$)}\returns true of $(yx)x^\rho=y$,}
\myline{\texttt{HasInverseProperty($L$)}\returns true if $x^\lambda(xy)=y=(yx)x^\rho$,}
\myline{\texttt{HasTwosidedInverses($L$)}\returns true if $x^\lambda=x^\rho$,}
\myline{\texttt{HasWeakInverseProperty($L$)}\returns true if $(xy)^\lambda x = y^\lambda$,}
\myline{\texttt{HasAutomorphicInverseProperty($L$)}\returns true if $(xy)^\lambda = x^\lambda y^\lambda$,}
\myline{\texttt{HasAntiautomorphicInverseProperty($L$)}\returns true if $(xy)^\lambda = y^\lambda x^\lambda$.}

\noindent Some properties of quasigroups:
\myline{\texttt{IsSemisymmetric($Q$)}\returns true if $(xy)x=y$,}
\myline{\texttt{IsTotallySymmetric($Q$)}\returns true if $Q$ is semisymmetric and commutative,}
\myline{\texttt{IsIdempotent($Q$)}\returns true it $x^2=x$,}
\myline{\texttt{IsSteinerQuasigroup($Q$)}\returns true if $Q$ is totally symm. and commutative,}
\myline{\texttt{IsUnipotent($Q$)}\returns true if $x^2=y^2$,}
\myline{\texttt{IsLeftDistributive($Q$)}\returns true if $x(yz) = (xy)(xz)$,}
\myline{\texttt{IsRightDistributive($Q$)}\returns true if $(xy)z=(xz)(yz)$,}
\myline{\texttt{IsDistributive($Q$)}\returns true if $Q$ is left and right distributive,}
\myline{\texttt{IsEntropic($Q$)}, \texttt{IsMedial($Q$)}\returns true if $(xy)(uv) = (xu)(yv)$.}

\noindent Loops of Bol-Moufang type:
\myline{\texttt{IsExtraLoop($L$)}\returns true if $x(y(zx)) = ((xy)z)x$,}
\myline{\texttt{IsCLoop($L$)}\returns true if $x(y(yz)) = ((xy)y)z$,}
\myline{\texttt{IsMoufangLoop($L$)}\returns true if $(xy)(zx)=(x(yz))x$,}
\myline{\texttt{IsRCLoop($L$)}\returns true if $x((yz)z) = (xy)(zz)$,}
\myline{\texttt{IsLCLoop($L$)}\returns true if $(xx)(yz) = (x(xy))z$,}
\myline{\texttt{IsRightBolLoop($L$)}\returns true if $x((yz)y) = ((xy)z)y$,}
\myline{\texttt{IsLeftBolLoop($L$)}\returns true if $x(y(xz)) = (x(yx))z$,}
\myline{\texttt{IsFlexible($Q$)}\returns true if $x(yx)=(xy)x$,}
\myline{\texttt{IsRightNuclearSquareLoop($L$)}\returns true if $x(y(zz)) = (xy)(zz)$,}
\myline{\texttt{IsMiddleNuclearSquareLoop($L$)}\returns true if $x((yy)z) = (x(yy))z$,}
\myline{\texttt{IsLeftNuclearSquareLoop($L$)}\returns true if $(xx)(yz) = ((xx)y)z$,}
\myline{\texttt{IsRightAlternative($Q$)}\returns true if $x(yy) = (xy)y$,}
\myline{\texttt{IsLeftAlternative($Q$)}\returns true if $(xx)y = x(xy)$,}
\myline{\texttt{IsAlternative($Q$)}\returns true if it is both left and right alternative.}

\noindent Power alternative loops:
\myline{\texttt{IsLeftPowerAlternative($L$)}\returns true if $x^n(x^my) = x^{n+m}y$,}
\myline{\texttt{IsRightPowerAlternative($L$)}\returns true if $(xy^n)y^m = xy^{n+m}$,}
\myline{\texttt{IsPowerAlternative($L$)}\returns true if $L$ is left and right power alternative.}

\noindent Conjugacy closed loops:
\myline{\texttt{IsLCCLoop($L$)}\returns true if left translations are closed under conjugation,}
\myline{\texttt{IsRCCLoop($L$)}\returns true if right translations are closed under conjugations,}
\myline{\texttt{IsCCLoop($L$)}\returns true if $L$ is left and right conjugacy closed.}

\noindent Additional varieties of loops:
\myline{\texttt{IsLeftBruckLoop($L$)}, \texttt{IsLeftKLoop($L$)}\returns true if $L$ is left Bol and has}
\myextraline{the automorphic inverse property,}
\myline{\texttt{IsRightBruckLoop($L$)}, \texttt{IsRightKLoop($L$)}\returns true if $L$ is right Bol and}
\myextraline{has the automorphic inverse property.}

\noindent Here is a nice, albeit trivial illustration of the filters built into the LOOPS package:

\begin{GAPexample}\begin{verbatim}
gap> L := LoopByCayleyTable([[1,2],[2,1]]);
<loop of order 2>
gap> IsLeftBolLoop(L); L;
true
<left Bol loop of order 2>
gap> IsRightBolLoop(L); L;
true
<Moufang loop of order 2>
gap> IsAssociative(L); L;
true
<associative loop of order 2>
\end{verbatim}\end{GAPexample}

\subsection{Libraries}

Several libraries of small loops up to isomorphism are included in LOOPS. As of
version 1.4.0, the libraries contain:
\begin{enumerate}
\item[-] all nonassociative left Bol loops of order $\le 16$,
\item[-] all nonassociative Moufang loops of order $\le 64$ and $= 81$,
\item[-] all nonassociative Steiner loops of order $\le 16$,
\item[-] all (three) nonassociative conjugacy closed loops of order $p^2$, for every odd prime $p$,
\item[-] all (one) nonassociative conjugacy closed loops of order $2p$, for every odd prime $p$,
\item[-] the smallest nonassociative simple Moufang loop (of order $120$),
\item[-] all nonassociative loops of order $\le 6$.
\end{enumerate}
There is also a library of all nonassociative loops of order $\le 6$ up to
isotopism.

The $m$th loop of order $n$ in a given library can be retrieved via
\begin{displaymath}
    \texttt{LeftBolLoop($n$,$m$)},\quad \texttt{MoufangLoop($n$,$m$)},
\end{displaymath}
and so on.

We took great care to store the information in the libraries efficiently. For
instance, the library of Moufang loops can be packed into less than $18$
kilobytes, hence averaging about $4$ bytes per loop.

\begin{remark} All nonassociative Moufang loops of order less than $64$ can be
found in \cite{Goodaire}. Our numbering for these loops agrees with
\cite{Goodaire}.

The 4262 nonassociative Moufang loops of order $64$ were first constructed in
\cite{PVEJC}, but it was proved (computationally) only in \cite{NagyVojt} that
the list is complete.

The 2038 nonassociative left Bol loops of order $16$ were enumerated for the
first time by Moorhouse \cite{Moorhouse}. The first author obtained the same
result by a different method, on which he will report in a separate paper
\cite{NagyBol}.

The fact that for every odd prime $p$ there are precisely three nonassociative
conjugacy closed loops of order $p^2$ was established by Kunen \cite{Kunen}.
Dr\'apal and Cs\"{o}rg\H{o} derived simple formulas for multiplication in
these three loops \cite{DrapalCsorgo}. When $p$ is an odd prime, Wilson
\cite{Wilson} constructed a nonassociative conjugacy closed loop of order $2p$,
and Kunen \cite{Kunen} showed there are no other such loops.

Our counts of small loops agree with the known results, e.g. \cite{McKay}.

The library of small Steiner loops is based on \cite{Rosa}.
\end{remark}

\section{Constructing isomorphisms}

There does not appear to be much research on the problem of finding an
isomorphism between loops. In this section we explain the approach used in
LOOPS. It works surprisingly well for many varieties of loops, including
Moufang loops.

Let $Q$ be a loop, and let $\mathcal P$ be a set of properties (of elements)
invariant under isomorphisms. The nature of $\mathcal P$ depends on $Q$. For
instance, when $Q$ is power-associative, one of the invariant properties for an
element $x$ might be the order $|x|$.

Given $\mathcal P$ and a collection $\mathcal C$ of loops, define an
equivalence on the (disjoint) union of $\mathcal C$ by $x\sim y$ if and only if
$\varphi(x)=\varphi(y)$ for every $\varphi\in\mathcal P$. Then, if $f:Q\to L$
is an isomorphism and $\mathcal C=\{Q$, $L\}$, we must have $x\sim f(x)$ for
every $x\in Q$. In other words, $\mathcal P$ partitions the elements into blocks
invariant under isomorphism.

\emph{In order to find an isomorphism, we need a set of invariants $\mathcal
P$ that is easy to calculate but results in a fine partition.}

We have used the following invariants $\mathcal P$ for power-associative loops:
\begin{align*}
    \varphi_1(x)&=|x|,\\
    \varphi_2(x)& = |\{y;\;y^2=x\}|,\\
    \varphi_3(x)& = |\{y;\;y^4=x\}|,\\
    \varphi_{4,k}(x) &= |\{y;\;xy=yx,\,|y|=k\}|,\text{ for $k\ge 1$}.
\end{align*}

The algorithm searching for an isomorphism $f:Q\to L$ first orders the
equivalence classes of $\sim$ by increasing size on both $Q$ and $L$. If the
equivalence class sizes of $Q$ and $L$ do not match, it is clear that no
isomorphism $f:Q\to L$ exists, and \texttt{fail} is returned. Otherwise, a
backtrack search attempts to find an isomorphism respecting the partitions of
$\sim$.

It would be an interesting project to analyze the speed of the algorithm
depending on the choice of $\mathcal P$. We do not claim that the above
$\mathcal P$ is optimized in any sense. Note, for instance, that the invariants
$\varphi_2$, $\varphi_3$ are useless for many power associative loops of odd order,
and $\varphi_{4,k}$ are useless for all commutative loops.

\section{Classification of small Frattini Moufang loops of order $64$}

Let $L$ be a loop and let the \emph{Frattini subloop} $\Phi(L)$ be the normal
subloop generated by all squares, commutators and associators of $L$. In other
words, $\Phi(L)$ is the smallest normal subloop such that $L/\Phi(L)$ is an
elementary abelian $p$-group. Following Hsu \cite{Hsu}, we say that $L$ is a
\emph{small Frattini $p$-loop} if $|\Phi(L)|\leq p$.

In this section, $L$ will denote a small Frattini Moufang 2-loop of order
$2^{n+1}$. Moreover, in order to avoid trivialities, we assume that
$|\Phi(L)|=2$. Clearly, $\Phi(L)\leq Z(L)$, $L$ is nilpotent of class $2$, and
it has a unique nontrivial square, commutator and associator element.

\begin{remark}
Small Frattini Moufang 2-loops are also called \emph{code loops} due to their
connection to doubly even linear binary codes. Some of these loops play an
important role in the description of large sporadic simple groups.
\end{remark}

We consider $V=L/\Phi(L)$ as a vector space over $\mathbb{F}_2$, and we
identify $\Phi(L)$ and $\mathbb{F}_2$. In particular, we sometimes write the
group operations additively.

Let us take arbitrary elements $u=x \mod \Phi(L)$, $v=y \mod \Phi(L)$, $w=z
\mod \Phi(L)$ of $V$. Then, the following maps are well defined:
\[ \begin{array}{ll}
    \sigma:V\to\mathbb{F}_2, & \sigma(u) = x^2, \\
    \kappa:V\times V\to\mathbb{F}_2, & \kappa(u,v) = [x,y], \\
    \alpha:V\times V\times V \to\mathbb{F}_2, & \alpha(u,v,w) = [x,y,z].
\end{array} \]
Moreover, $\alpha$ is an alternating trilinear form, $\kappa$ is alternating,
and we have
\begin{align*}
    &\sigma(u+v)=\sigma(u)+\sigma(v)+\kappa(u,v),\\
    &\kappa(u+v,w)=\kappa(u,w)+\kappa(v,w)+\alpha(u,v,w).
\end{align*}
Hence, by definition, $V$ is a \emph{symplectic cubic space.}

There are different ways in which a small Frattini Moufang $2$-loop is obtained
from a symplectic cubic space (cf. Griess \cite{Griess}, Chein and Goodaire
\cite{CheinGoodaire}, Hsu \cite{Hsu}). All of the above constructions induct on
the dimension of $V$. In contrast, a new approach, \cite{NagyCL}, takes
advantage of \emph{groups with triality} and constructs the loop globally.

For this, let $\sigma_i$, $\kappa_{ij}$ and $\alpha_{ijk}$ be the structure
constants of $\sigma, \kappa, \alpha$ with respect to a fixed basis ov $V$. We
define the group $G$ with gerenators $g_i, f_i, h_i$, $i\in \{1,\ldots, n\}$,
$u$ and $v$ by the following relations:
\begin{eqnarray*}
&& g_i^2 = u^{\sigma_i}, \; f_i^2 = v^{\sigma_i}, \; h_i^2=u^2=v^2=1,
  \label{eq:G1} \\
&& [g_i,g_j] = u^{\kappa_{ij}}, \; [f_i,f_j] = v^{\kappa_{ij}},
  \label{eq:G2} \\
&& [g_i,f_j] = (uv)^{\kappa_{ij}} \, \prod_{k=1}^n h_k^{\alpha_{ijk}},
  \label{eq:G3} \\
&& [g_i,h_j] = u^{\delta_{ij}}, \; [f_i,h_j] = v^{\delta_{ij}},
  \label{eq:G4} \\
&& [h_i,h_j] = [g_i,u] = [f_i,u] = [h_i,u] = [g_i,v] = [f_i,v] = [h_i,v] =
  1. \label{eq:G5}
\end{eqnarray*}
Then, $G$ is a group and the maps
\begin{eqnarray*}
\tau&:& g_i \leftrightarrow f_i, \, h_i \mapsto h_i, \, u \leftrightarrow
v \label{eq:sigma} \\
\rho&:& g_i \mapsto f_i, \, f_i \mapsto (g_if_i)^{-1}, \, h_i \mapsto h_i,
\, u \mapsto v, \, v \mapsto uv \label{eq:rho}
\end{eqnarray*}
extend to \emph{triality automorphisms} of $G$. The following function returns
the Moufang loop associated to the group $G$ with triality automorphisms $\tau,
\rho$:
\begin{GAPexample}\begin{verbatim}
TrialityGroupToLoop := function( G, tau, rho )
   local ccl, ct;
   ccl := Elements( ConjugacyClass( G, tau ) );
   ct := List( ccl, s1 ->
            List( ccl, s2 ->
               Position( ccl, s1^rho * s2^(rho^2) * s1^rho )
            )
         );
   return LoopByCayleyTable( NormalizedQuasigroupTable( ct ) );
end;
\end{verbatim}\end{GAPexample}
To complete the classification of small Frattini Moufang 2-loops of order $64$,
it now suffices to classify the symplectic cubic spaces of order $32$. For a
fixed basis, such a space is given by
\begin{displaymath}
    \left ( \begin{array}{c} 5\\3 \end{array} \right ) +
    \left ( \begin{array}{c} 5\\2 \end{array} \right ) + 5 =25
\end{displaymath}
structure constants, which give rise to a 25-dimensional vector space $W$ over
$\mathbb{F}_2$.

Any linear map $A$ of $V$ defines a new symplectic cubic space with maps
\[\sigma^A(u)=\sigma(Au),\quad \kappa^A(u,v)=\kappa(Au,Av),\quad \alpha^A(u)=\alpha(Au,Av,Aw),\]
and hence $A$ induces a linear map on $W$. This defines an action of $GL(5,2)$
on $W$.

It is easy to show the 1-1 correspondence of loop isomorphisms and linear
isomorphisms of symplectic cubic spaces. This implies that the orbits of
$GL(5,2)$ on $W$ will correspond precisely to the isomorphism classes of small
Frattini Moufang 2-loops of order $64$.

Since $|GL(5,2)|$ and $2^{25}$ are still too large for GAP to compute the needed
orbits, one has to have a closer look at invariant subspaces of $W$. Once this
is done, the classification is complete, with the result that there are
precisely $80$ nonisomorphic small Frattini Moufang loops of order $64$.

\section{An interesting Cs\"{o}rg\H{o} loop}

One of the longer-standing problems in loop theory was the question if there is
a loop with nilpotency class higher than two whose inner mapping group is
abelian. In \cite{Csorgo}, Cs\"{o}rg\H{o} constructed such a loop (of order
$128$ and nilpotency class $3$). The following GAP code returns this loop $L$.
The code follows \cite{Csorgo}, where some insight is given.

\begin{GAPexample}\begin{verbatim}
# constructing a group of order 8192 by presenting relations
f := FreeGroup(13);
G := f/[ f.1^2, f.2^2, f.3^2, f.4^2, f.5^2, f.6^2, f.7^2, f.8^2, f.9^2, f.10^2,
f.11^2, f.12^2, f.13^2, (f.1*f.2)^2, (f.1*f.3)^2, (f.1*f.4)^2, (f.1*f.5)^2,
(f.1*f.6)^2, (f.1*f.7)^2, (f.1*f.8)^2, (f.1*f.9)^2, (f.1*f.10)^2, (f.1*f.11)^2,
(f.1*f.12)^2, (f.1*f.13)^2, (f.2*f.3)^2, (f.2*f.4)^2, (f.3*f.4)^2, (f.2*f.5)^2,
(f.2*f.6)^2, (f.2*f.7)^2, (f.3*f.5)^2, (f.3*f.6)^2, (f.3*f.7)^2, (f.4*f.5)^2,
(f.4*f.6)^2, (f.4*f.7)^2, (f.2*f.9)^2, (f.2*f.10)^2, (f.3*f.8)^2, (f.3*f.10)^2,
(f.4*f.8)^2, (f.4*f.9)^2, f.1*f.2*f.8*f.2*f.8, f.1*f.3*f.9*f.3*f.9,
f.1*f.4*f.10*f.4*f.10, (f.2*f.11)^2, (f.2*f.12)^2, (f.2*f.13)^2, (f.3*f.11)^2,
(f.3*f.12)^2, (f.3*f.13)^2, (f.4*f.11)^2, (f.4*f.12)^2, (f.4*f.13)^2, (f.5*f.6)^2,
(f.5*f.7)^2, (f.6*f.7)^2, (f.5*f.9)^2, (f.5*f.10)^2, (f.6*f.8)^2, (f.6*f.10)^2,
(f.7*f.8)^2, (f.7*f.9)^2, f.1*f.5*f.8*f.5*f.8, f.1*f.6*f.9*f.6*f.9,
f.1*f.7*f.10*f.7*f.10, (f.5*f.12)^2, (f.5*f.13)^2, (f.6*f.11)^2, (f.6*f.13)^2,
(f.7*f.11)^2, (f.7*f.12)^2, f.1*f.11*f.5*f.11*f.5, f.1*f.12*f.6*f.12*f.6,
f.1*f.13*f.7*f.13*f.7, f.2*f.5*f.9*f.10*f.9*f.10, f.3*f.6*f.8*f.10*f.8*f.10,
f.4*f.7*f.8*f.9*f.8*f.9, (f.8*f.11)^2, (f.9*f.12)^2, (f.10*f.13)^2,
f.8*f.12*f.8*f.4*f.12*f.7, f.8*f.13*f.8*f.3*f.13*f.6, f.10*f.11*f.10*f.3*f.11*f.6,
f.9*f.11*f.9*f.11*f.7, f.9*f.13*f.9*f.13*f.5, f.10*f.12*f.10*f.12*f.5,
(f.11*f.12)^2, (f.11*f.13)^2, (f.12*f.13)^2 ];
# auxiliary data
g := GeneratorsOfGroup(G);
N := Subgroup( G, [ g[5], g[6], g[7], g[1] ] );
W := Subgroup( G, [ g[5]*g[2], g[6]*g[3], g[7]*g[4], g[1] ] );
A_0 := [ One(G), g[8], g[9], g[10], g[8]*g[9], g[8]*g[10], g[9]*g[10]*g[2],
    g[8]*g[9]*g[10]*g[2] ];
B_0 := [ One(G), g[8]*g[11], g[9]*g[12], g[10]*g[13], g[8]*g[11]*g[9]*g[12],
    g[8]*g[11]*g[10]*g[13]*g[3], g[9]*g[12]*g[10]*g[13],
    g[8]*g[11]*g[9]*g[12]*g[10]*g[13]*g[3] ];
A := Union( List( Elements( N ), x -> A_0*x ) );
B := Union( List( Elements( W ), x -> B_0*x ) );
H := Subgroup( G, [ g[2], g[3], g[4], g[11], g[12], g[13] ] );
# constructing the loop
ListPosition := function( S, x )
    local i; i := 1; while not x in S[i] do i := i + 1; od; return i;
end;
m := MappingByFunction( Domain(Elements( G)), Domain([1..8192]),
    x -> Position( Elements(G), x ) );
CA := List( A, x -> x*Elements( H ) );
mCA := List( CA, c -> Set( c, x -> x^m ) );
T := List([1..128],i->[1..128]);
for ii in [1..128] do for jj in [1..128] do
    T[ii][jj] := ListPosition( mCA, (A[ii]*B[jj])^m );
od; od;
p := SortingPerm( T[1] );
T := List( T, r -> Permuted( r, p ) );
L := LoopByCayleyTable( T );
\end{verbatim}\end{GAPexample}

In addition, the following properties hold for $L$: (a) the nucleus of $L$ is
elementary abelian of order $16$, (b) the left and middle nuclei have order
$32$, (c) the right nucleus has order $16$, (d) the two-element center
coincides with the associator subloop.

An interesting, more symmetric loop $K$ is obtained from $L$ by this greedy
algorithm:

Given a groupoid $Q$, let $\mu(Q) = |\{(a,b,c)\in Q\times Q\times Q;\;a(bc)\ne
(ab)c\}|$. Hence $\mu(Q)$ is a crude measure of (non)associativity of $Q$.

Let $T$ be a multiplication table of $L$ split into blocks of size $16\times
16$ according to the cosets of the nucleus of $L$. Let $h$ be the nontrivial
central element of $L$.

(*) For $1\le i\ne j\le 16$, let $T_{ij}$ be obtained from $T$ by multiplying
the $(i,j)$th block and the $(j,i)$th block of $T$ by $h$. Let $(s,t)$ be such
that $\mu(T_{st})$ is minimal among all $\mu(T_{ij})$. If $\mu(T_{st})
\ge\mu(T)$, stop, and return $T$. Else replace $T$ by $T_{st}$, and repeat (*).

It turns out that the multiplication table $T$ found by the above greedy
algorithm yields another loop $K$ of nilpotency class $3$ whose inner mapping
group is abelian. In addition, the following properties hold for $K$: (a) the
nucleus is elementary abelian of order $16$, (b) the left, middle, and right
nuclei have order $64$, (c) the two-element center coincides with the
associator subloop. In particular, $K$ is not isomorphic to $L$. Among other
peculiar features, it contains a nonassociative power associative loop of order
$64$ that is the union of its nuclei.

The construction of $L$ takes a minute or so in GAP, since calculations in free
groups are slow. A more direct, systematic, and much faster construction of $L$
and $K$ will be presented elsewhere \cite{DV}.

\noindent \footnotesize{\hfill Received \ February 25, 2007

\medskip
\begin{tabular}{lll}
        G\'abor P. Nagy&\quad&Petr Vojt\v{e}chovsk\'y\\
        Bolyai Institute&\quad&Department of Mathematics\\
        University of Szeged&\quad&University of Denver\\
        Aradi v\'ertan\'uk tere 1&\quad&2360 S Gaylord St\\
        H-6720 Sze\-ged&\quad&Denver, Colorado 80208\\
        Hungary&\quad&U.S.A.\\
        e-mail: nagyg@math.u-szeged.hu&\quad&e-mail: petr@math.du.edu
\end{tabular}
}

\end{document}